\RequirePackage[english]{babel} 

\documentclass{etds}
\usepackage[utf8]{inputenc}      
\usepackage[T1]{fontenc}         
\usepackage[varg]{txfonts}
\usepackage[dvipsnames]{xcolor}  
\usepackage{tikz}
\usepackage{Z-tikz}

\theoremstyle{plain}
\newtheorem{theo}{Theorem}
\newtheorem{lemm}{Lemma}

\usepackage[backend=biber,bibstyle=Z-bib,citestyle=alphabetic]{biblatex}

\usepackage[colorlinks=true,%
            linkcolor=RoyalBlue,%
            citecolor=RedViolet,%
            urlcolor=RoyalBlue]{hyperref}   

  \newcommand\CC{{\mathbf{C}}}         
  \newcommand\RR{{\mathbf{R}}}         
  \newcommand\TT{{\mathbf{T}}}         
  \newcommand\ZZ{{\mathbf{Z}}}         
  \newcommand\Lie{\mathfrak}
  \newcommand\haar{\eta}
  \newcommand\kprime{k'}
\renewcommand\O{{\mathcal O}}
  \newcommand\poids\nu
  \newcommand{\<}{\langle}
\renewcommand{\>}{\rangle}

\bibliography{biblio}

\begin{document}
\ETDS{1}{5}{00}{0000}
\runningheads{R. Howe and F. Ziegler}{Bohr density of simple linear group orbits}

\title{Bohr density of simple linear group orbits}

\author{ROGER HOWE\affil{1}\ and FRANÇOIS ZIEGLER\affil{2}}
\address{
\affilnum{1}\ 
Department of Mathematics, Yale University, New Haven, CT 06520-8283, USA\\
\email{howe@math.yale.edu}\\
\affilnum{2}\ 
Department of Mathematical Sciences, Georgia Southern University, Statesboro,\\ GA 30460-8093, USA\\
\email{fziegler@georgiasouthern.edu}}

\recd{$15$ November $2012$ and accepted in revised form $26$ July $2013$}


\begin{abstract}
   We show that any nonzero orbit under a noncompact, \mbox{simple}, irreducible linear group is dense in the Bohr compactification of the ambient space.
\end{abstract}

\section{Introduction}
Let $V$ be a locally compact abelian group, $V^*$ its Pontrjagin dual and  $bV$ its Bohr compactification, i.e.~$bV$ is the dual of the discretized group $V^*$. On identifying $V$ with its double dual we have a dense embedding $V\hookrightarrow bV$, viz.
\begin{equation*}
   \{\textrm{continuous characters of $V^*$}\}
   \hookrightarrow
   \{\textrm{all characters of $V^*$}\}.
\end{equation*}
The relative topology of $V$ in $bV$ is known as the \emph{Bohr topology} of $V$. Among its many intriguing properties (surveyed in \cite{Galindo:2007}) is the observation due to Katznelson \cites{Katznelson:1973}[§7.6]{Graham:1979} that very ``thin'' subsets of $V$ can be Bohr dense in very large ones. 

While Katznelson was concerned with the case $V=\ZZ$ (the integers), we shall illustrate this phenomenon in the setting where $V$ is the additive group of a real vector space, and the subsets of interest are the orbits of a Lie group acting linearly on $V$. Indeed our aim is to establish the following result, which was conjectured in \cite[p.\,45]{Ziegler:1996}:

\begin{theo}
   \label{SimpleThm}
   Let $G$ be a noncompact\textup, simple real Lie group and $V$ a nontrivial\textup, irreducible\textup, finite-di\-men\-sional real $G$-module. Then every nonzero $G$-orbit in $V$ is dense in $bV$. 
\end{theo}

\noindent
We prove this in §3 on the basis of four lemmas prepared in §2. Before that, let us record a similar property of \emph{nilpotent} groups. In that case, orbits typically lie in proper affine subspaces, so we can't hope for Bohr density in the whole space; but we have:

\begin{theo}
   \label{NilpotentThm}
   Let $G$ be a connected nilpotent Lie group and $V$ a finite-dimen\-sional $G$-module of unipotent type. Then every $G$-orbit in $V$ is Bohr dense in its affine hull. 
\end{theo}

\proc{Proof.}
   Recall that \emph{unipotent type} means that the Lie algebra $\Lie g$ of $G$ acts by nilpotent operators. So $Z\mapsto\exp(Z)v$ is a polynomial map of $\Lie g$ onto the orbit of $v\in V$, and the claim follows immediately from \cite[Theorem]{Ziegler:1993}. \ep

\section{Four lemmas}

Our first lemma gives several characterizations of Bohr density --- each of which can also be regarded as providing a corollary of Theorem \ref{SimpleThm}.

\begin{lemm}
   \label{Lemma1}
   Let $\O$ be a subset of the locally compact abelian group $V$. Then the following are equivalent\textup:
   \begin{enumerate}
      \item[$(1)$] $\O$ is dense in $bV$\textup;
      \item[$(2)$] $\alpha(\O)$ is dense in $\alpha(V)$ whenever $\alpha$ is a continuous morphism from $V$ to a compact topological group\textup;
      \item[$(3)$] Every almost periodic function on $V$ is determined by its restriction to~$\O$\textup;
      \item[$(4)$] Haar measure $\haar$ on $bV$ is the weak$*$ limit of probability measures $\mu_T$ concentrated on $\O$.
   \end{enumerate}
\end{lemm}

\proc{Proof.}
   $(1)\Leftrightarrow(2)$: Clearly (2) implies (1) as the special case where $\alpha$ is the natural inclusion $\iota:V\hookrightarrow bV$. Conversely, suppose (1) holds and $\alpha:V\to X$ is a continuous morphism to a compact group. By the universal property of $bV$ \cite[{}16.1.1]{Dixmier:1982}, $\alpha=\beta\circ\iota$ for a continuous morphism $\beta:bV\to X$. Now continuity of $\beta$ implies $\beta(\overline{\iota(\O)})\subset\overline{\beta(\iota(\O))}$, which is to say that $\beta(bV)\subset\overline{\alpha(\O)}$ and hence $\alpha(V)\subset\overline{\alpha(\O)}$, as claimed.

   $(1)\Leftrightarrow(3)$: Recall that a function on $V$ is \emph{almost periodic} iff it is the pull-back of a continuous $f:bV\to\CC$ by the inclusion $V\hookrightarrow bV$. If two such functions coincide on $\O$ and $\O$ is dense in $bV$, then clearly they coincide everywhere. Conversely, suppose that $\O$ is not dense in $bV$. Then by complete regularity \cite[{}8.4]{Hewitt:1963} there is a nonzero continuous $f:bV\to\CC$ which is zero on the closure of $\O$ in $bV$. Now clearly this $f$ is not determined by its restriction to $\O$.
  
   $(1)\Leftrightarrow(4)$ (\cite{Katznelson:1973}): Suppose that $\haar$ is the weak$*$ limit of probability measures $\mu_T$ concentrated on $\O$. So we have $\mu_T(f)\to\haar(f)$ for every continuous $f$, and the complement of $\O$ in $bV$ is $\mu_T$-null \cites[Def.~V.5.7.4 and Prop.~IV.5.2.5]{Bourbaki:2004}. If $f$ vanishes on the closure of $\O$ in $bV$ then so do all $\mu_T(|f|)$ and hence also $\haar(|f|)$, which forces $f$ to vanish everywhere. So $\O$ is dense in $bV$. Conversely, suppose that $\O$ is dense in $bV$. We have to show that given continuous functions $f_1,\dots,f_n$ on $bV$ and $\varepsilon>0$, there is a probability measure $\mu$ concentrated on $\O$ such that $|\haar(f_j)-\mu(f_j)|<\varepsilon$ for all $j$. Writing
   \begin{equation*}
      F=(f_1,\dots,f_n)
      \quad\text{and}\quad
      \haar(F)=(\haar(f_1),\dots,\haar(f_n))
   \end{equation*}
   we see that this amounts to $\|\haar(F)-\mu(F)\|<\varepsilon$, where the norm is the sup norm in $\CC^n$. Now by \cite[Cor.~V.6.1]{Bourbaki:2004} $\haar(F)$ lies in the convex hull of $F(bV)$ (which is compact by Carathéodory's theorem \cite[{}11.1.8.7]{Berger:1987}). So $\haar(F)$ is a convex combination $\smash{\sum_{i=1}^N\lambda_iF(\omega_i)}$ of elements of $F(bV)$. But $F(\O)$ is dense in $F(bV)$, so we can find $w_i\in\O$ such that $\|F(\omega_i)-F(w_i)\|<\varepsilon$. Putting $\mu=\smash{\sum_{i=1}^N\lambda_i\delta_{w_i}}$ where $\delta_{w_i}$ is Dirac measure at $w_i$, we obtain the desired probability measure $\mu$. \ep
\medbreak

\proc{Remark 1.}
   One might wonder if condition (2) is equivalent to the following \emph{a priori} weaker but already interesting property:
   \begin{enumerate}
      \item[$(2')$] $\O$ has dense image in any compact quotient group of $V$.
   \end{enumerate}
   Here is an example showing that $(2')$ \emph{does not} imply (2): Let $V=\RR$ and $\O = \ZZ\cup2\pi\ZZ$. Then clearly $\O$ has dense image in every compact quotient $\RR/a\ZZ$. On the other hand, considering the irrational winding $\alpha:\RR\to\TT^2$ defined by $\alpha(v)=(\mathrm{e}^{\mathrm{i}v},\mathrm{e}^{2\pi\mathrm{i}v})$, one checks without trouble that  $\overline{\alpha(\O)}=\TT\times\{1\}\cup\{1\}\times\TT$, which is strictly smaller than $\overline{\alpha(V)}=\TT^2$.
\medbreak

\proc{Remark 2.}
   A net of probability measures $\mu_T$ converging to Haar measure on $bV$ as in (4) has been called a \emph{generalized summing sequence} by Blum and Eisenberg \cite{Blum:1974}. They observed, among others, the following characterization.
\medbreak

\begin{lemm}
   \label{Lemma2}
   The following conditions are equivalent\textup:
   \begin{enumerate}
      \item[$(1)$] $\mu_T$ is a generalized summing sequence\textup;
      \item[$(2)$] The Fourier transforms $\hat\mu_T(u)=\int_{bV}\omega(u)d\mu_T(\omega)$ converge pointwise to the characteristic function of $\smash{\{0\}\subset V^*}$.
   \end{enumerate}
\end{lemm}

\proc{Proof.}
   This characteristic function is the Fourier transform of Haar measure $\haar$ on $bV$. Thus, condition (2) says that $\mu_T(f)\to\haar(f)$ for every continuous character $f(\omega)=\omega(u)$ of $bV$; whereas condition (1) says that $\mu_T(f)\to\haar(f)$ holds for every continuous function $f$ on $bV$. Since linear combinations of continuous characters are uniformly dense in the continuous functions on $bV$ (Stone-Weierstrass), the two conditions imply each other. \ep
\medbreak

For our third lemma, let $G$ be a group, $V$ a finite-dimensional $G$-mod\-ule, and write $V^*$ for the dual module wherein $G$ acts contragrediently: $\<gu,v\>=\<u,g^{-1}v\>$. We have

\begin{lemm}
   \label{Lemma3}
   Suppose that $V$ is irreducible and $\phi(g)=\<u,gv\>$ is a nonzero matrix coefficient of $V$. Then every other matrix coefficient $\psi(g)=\<x,gy\>$ is a linear combination of left and right translates of $\phi$.
\end{lemm}

\proc{Proof.}
   Irreducibility of $V$ and (therefore) $V^*$ ensures that $u$ and $v$ are cyclic, i.e.~their $G$-orbits span $V^*$ and $V$. So we can write $x=\smash{\sum_i \alpha_ig_iu}$ and $y=\smash{\sum_j \beta_jg_jv}$, whence $\psi(g)=\smash{\sum_{i,j} \alpha_i\beta_j\phi(g_{\smash{\raisebox{1pt}{$\scriptstyle i$}}}^{-1}gg_j)}$. \ep
\medbreak 

Finally, our fourth preliminary result is the famous

\begin{lemm}[Van der Corput]
   \label{LemmaVdC}
   Suppose that $F:[a,b]\to\RR$ is differentiable\textup, its derivative $F'$ is monotone\textup, and $|F'|\geqslant 1$ on $(a,b)$. Then $\bigl|\int_a^b\mathrm{e}^{\mathrm{i}F(t)}dt\,\bigr|\leqslant3$.
\end{lemm}

\proc{Proof.}
   See \cite[p.\,332]{Stein:1993}, or \cite[Lemma 3]{Rogers:2005} which actually gives the sharp bound 2. \ep
\medbreak

\section{Proof of Theorem \ref{SimpleThm}}
By Lemma \ref{Lemma1}, it is enough to show that Haar measure on $bV$ is the weak$*$ limit of probability measures $\mu_T$ concentrated on the orbit under consideration; or equivalently (Lemma \ref{Lemma2}), that the Fourier transforms of the $\mu_T$ tend pointwise to the characteristic function of $\{0\}\subset V^*$. (Here we identify the Pontrjagin dual with the dual vector space or module.)
   
To construct such $\mu_T$, we assume without loss of generality that the action of $G$ on $V$ is effective, so that we may regard $G\subset\operatorname{GL}(V)$. Let $K\subset G$ be a maximal compact subgroup, $\Lie g = \Lie k + \Lie p$ a Cartan decomposition, $\Lie a\subset \Lie p$ a maximal abelian subalgebra, $C\subset\Lie a^*$ a Weyl chamber, $P\subset\Lie a$ the dual positive cone, and $H$ an interior point of $P$; thus we have that $\<\poids,H\>$ is positive for all nonzero $\poids\in C$. (For all this structure see, for example, \cite{Kostant:1973}.) We fix a nonzero $v\in V$, and for each positive $T\in\RR$ we let $\mu_T$ denote the image of the product measure Haar $\times$ (Lebesgue$/T$) $\times$ Haar under the composed map
\begin{equation*}
   \raisebox{-4.8ex}{
   \begin{tikzpicture}[semithick,>=to]
      \node[anchor=east] (11) at (-1.8,0.6) {$K \times [0,T] \times K$};
      \node              (12) at (0,0.6) {$Gv$};
      \node[anchor=west] (13) at (1.8,0.6) {$bV$};
      \node[anchor=east] (21) at (-1.8,0) {$(k,    t,    \kprime)$};
      \node              (22) at (0,0) {\raisebox{1.1mm}{$k\exp(tH)\kprime v$}};
      \node              (32) at (0.0,-0.6) {$w$};
      \node[anchor=west] (33) at (1.8,-0.6) {\raisebox{1.7mm}{$\mathrm{e}^{\mathrm{i}\<\cdot,w\>}.$}};
      \draw[->] (11) -- (12);
      \draw[->] (12) -- (13);
      \draw[|->] (21) -- (22);
      \draw[|->] (32) -- (33);
   \end{tikzpicture}
   }
\end{equation*}
Here $\exp:\Lie a\to A$ is the usual matrix exponential with inverse $\log:A\to\Lie a$, and the brackets $\<\cdot,\cdot\>$ denote both pairings, $\Lie a^*\times \Lie a\to\RR$ and $V^*\times V\to\RR$. By construction the $\mu_T$ are concentrated on the subset $Gv$ of $bV$ \cite[Cor.~V.6.2.3]{Bourbaki:2004}. There remains to show that as $T\to\infty$ we have, for every nonzero $u\in V^*$,
\begin{equation}
   \label{limit}
   \int_{K \times K} dk\,d\kprime \,\frac1T\int_0^T \mathrm{e}^{\mathrm{i}\<u,k\exp(tH)\kprime v\>} dt \to 0.
   \tag{$*$}
\end{equation}
To this end, let
\begin{equation*}
   F_{k\kprime }(t)=\<u,k\exp(tH)\kprime v\>
\end{equation*}
denote the exponent in \eqref{limit}. We are going to show that Lemma \ref{LemmaVdC} applies to almost every $F_{k\kprime}$. In fact, it is well known (see for example \cite[Prop. 2.4 and proof of Prop.\,3.4]{Kostant:1973}) that $\Lie a$ acts diagonalizably (over $\RR$) on $V$. Thus, letting $E_\poids$ be the projector of $V$ onto the weight $\poids$ eigenspace of $\Lie a$, we can write 
\begin{equation*}
   F_{k\kprime }(t) = \sum_{\poids\in \Lie a^*} \<u,kE_{\poids}\kprime v\> \mathrm{e}^{\<\poids,H\>t}.
\end{equation*}
   Now we claim that there are nonzero $\poids$ such that the coefficient  $f_\poids(k,\kprime ) = \<u,kE_{\poids}\kprime v\>$ is not identically zero on $K \times K$. (Then $f_\poids$, being analytic, will be nonzero \emph{almost everywhere}.) Indeed, suppose otherwise. Then, writing any $g\in G$ in the form  $ka\kprime $ ($KAK$ decomposition \cite{Knapp:2002}), we would have
\begin{equation*}
   \<u,gv\> = \sum_{\poids\in\Lie a^*} \<u,kE_{\poids}\kprime v\> \mathrm{e}^{\<\poids,\log(a)\>} = \<u,kE_{0}\kprime v\>.
\end{equation*}
In particular the matrix coefficient $\<u,gv\>$ would be bounded. Hence so would be all matrix coefficients, since they are linear combinations of translates of this one (Lemma \ref{Lemma3}); and this would contradict the noncompactness of $G\subset\operatorname{GL}(V)$. 
   
So the set $N=\left\{\poids\in\Lie a^*:\poids\ne0,f_\poids\ne0\right\}$ is not empty. It is also Weyl group invariant, hence contains weights $\poids\in C$ for which we know $\<\poids,H\>$ is positive. Therefore, maximizing $\<\poids,H\>$ over $N$ produces a positive number $\<\poids_0,H\>$, in terms of which our exponent and its derivatives can be written
\begin{equation*}
   \frac{d^n}{dt^n}F_{k\kprime }(t) = \mathrm{e}^{\<\poids_0,H\>t} 
   \sum_{\poids\in\Lie a^*}f_\poids(k,\kprime)\<\poids,H\>^n
   \mathrm{e}^{\<\poids-\poids_0,H\>t}
\end{equation*}
where $\<\poids-\poids_0,H\><0$ in all nonzero terms except the one indexed by $\poids_0$. (Here we assume, as we may, that $H$ was initially chosen outside the kernels of all pairwise differences of weights of $V$.) From this it is clear that for almost all $(k,\kprime )$ there is a $T_0$ beyond which the first two derivatives of $F_{k\kprime}$ are greater than 1 in absolute value. So Lemma~\ref{LemmaVdC} applies and gives
\begin{equation*}
   \biggl|\int_{T_0}^T \mathrm{e}^{\mathrm{i}F_{k\kprime }(t)} dt\,\biggr| \leqslant 3
   \qquad\forall\,T.
\end{equation*}
Therefore we have $\lim_{T\to\infty}\frac1T\int_0^T \mathrm{e}^{\mathrm{i}F_{k\kprime }(t)} dt = 0$ for almost all $(k,\kprime )$, whence the conclusion \eqref{limit} by dominated convergence. This completes the proof.

\section{Outlook}

Theorem \ref{SimpleThm} says that the $G$-action on $V\setminus\{0\}$ is \emph{minimal} \cite{Petersen:1983} in the Bohr topology. It would be interesting to determine if it is still minimal, and/or \emph{uniquely ergodic}, on $bV\setminus\{0\}$.

It is also natural to speculate whether our theorems have a common extension to more general group representations. Here we shall content ourselves with noting two obstructions. First, Theorem \ref{SimpleThm} clearly fails for \emph{semi}simple groups with compact factors. Secondly, Theorem \ref{NilpotentThm} fails for $V$ not of unipotent type, as one sees by observing that the orbits of $\RR$ acting on $\RR^2$ by 
$
\smash{\exp
   \left(
      \begin{smallmatrix}
         t&\phantom{-}0\\
         0&-t
      \end{smallmatrix}
   \right)}
$
(i.e., hyperbolas) already have non-dense images in $\RR^2/\ZZ^2$.

\ack{We thank Francis Jordan, who found the example in Remark~1.} 

\setlength{\labelalphawidth}{3.2em}

\printbibliography

\end{document}